\documentclass{article}%
\usepackage{amsmath}
\usepackage{amsfonts}
\usepackage{amssymb}
\usepackage{graphicx}
\usepackage{setspace}
\usepackage{hyperref}
\hypersetup{colorlinks=true, linkcolor=blue, citecolor=red, urlcolor=red}
\newtheorem{theorem}{Theorem}

\newtheorem{corollary}[theorem]{Corollary}

\newtheorem{definition}[theorem]{Definition}

\newtheorem{lemma}[theorem]{Lemma}

\newenvironment{proof}[1][Proof]{\noindent\textbf{#1.} }{\ \rule{0.5em}{0.5em}}

\begin{document}
\doublespacing
\title{ \textbf{On Metrizing Vague Convergence of Random Measures with Applications on Bayesian Nonparametric Models}}   
\author{Luai Al-Labadi\thanks{{\em Address for correspondence:} Luai Al-Labadi, Department of Mathematical \& Computational Sciences, University of Toronto Mississauga, 3359 Mississauga Road, Deerfield Hall, Room 3015, Mississauga, Ontario M5S 3G3, Canada. E-mail: luai.allabadi@utoronto.ca.}}

\date{\vspace{-5ex}}
\maketitle
\pagestyle {myheadings} \markboth {} {Metrizing Vague Convergence}
\begin{abstract}
  This paper deals with studying  vague convergence of random measures of the form  $\mu_{n}=\sum_{i=1}^{n} p_{i,n} \delta_{\theta_i}$,
 where $(\theta_i)_{1\le i \le n}$  is a sequence of independent and identically distributed random variables   with common distribution $\Pi$, $(p_{i,n})_{1 \le i \le n}$ are random variables  chosen according to certain procedures and are independent of $(\theta_i)_{i \geq 1}$ and $\delta_{\theta_i}$ denotes the Dirac measure at $\theta_i$.  We show that  $\mu_{n}$ converges vaguely to $\mu=\sum_{i=1}^{\infty} p_{i} \delta_{\theta_i}$ if and only if $\mu^{(k)}_{n}=\sum_{i=1}^{k} p_{i,n} \delta_{\theta_i}$ converges vaguely  to $\mu^{(k)}=\sum_{i=1}^{k} p_{i} \delta_{\theta_i}$ for all $k$ fixed. The limiting process $\mu$  plays a central role in many areas in statistics, including Bayesian nonparametric models.  A finite approximation of the beta process is derived from the application of this result. A simulated example is incorporated, in which the proposed approach exhibits an excellent performance over several existing algorithms.

\par
 \vspace{9pt} \noindent\textsc{Key words:}  Beta process, Nonparametric Bayesian statistics,  Point Processes, Random measures,  Vague convergence.

 \vspace{9pt}

\noindent { \textbf{MSC 2000}} Primary 62F15, 60G57; Secondary 28A33.

\end{abstract}

\section {Introduction}
\label{intro}

The primary objective of this paper is to study the vague convergence of a particular class of random measures of the form
 \begin{equation}
 \mu_{n}=\sum_{i=1}^{n} p_{i,n} \delta_{\theta_i},
 \end{equation}\label{Pn}
 where $(\theta_i)_{i \geq 1}$  is a sequence of independent and identically distributed (i.i.d.) random variables   with common distribution $\Pi$, $(p_{i,n})_{1 \le i \le n}$ are random variables  chosen according to certain procedures and are independent of $(\theta_i)_{i \geq 1}$, and $\delta_{\theta_i}$ denotes the Dirac measure at $\theta_i$. In particular, we show that  $\mu_{n}$ converges vaguely to $\mu=\sum_{i=1}^{\infty} p_{i} \delta_{\theta_i}$ if and only if $\mu^{(k)}_{n}=\sum_{i=1}^{k} p_{i,n} \delta_{\theta_i}$ converges vaguely to $\mu^{(k)}=\sum_{i=1}^{k} p_{i} \delta_{\theta_i}$ for all $k$ fixed. The limiting process $\mu$ has played a central role in some area in statistics, including Bayesian nonparametric models. Interesting examples of $\mu$ include, among others, the Dirichlet process (Ferguson, 1973), the beta process  (Hjort, 1990), the beta-Stacy  process (Walker and Muliere, 1997), the two parameter Poisson-Dirichlet process (Pitman and Yor, 1997) and the normalized inverse-Gaussian process (Lijoi, Mena and Pr\"unster, 2005). For a recent summary of the Bayesian nonparametric priors and their applications, please refer to the book of Phadia (2013) and M\"uller, Quintana, Jara and Hanson  (2015). A motivating  application of the aforementioned result involves the derivation of a finite sum representation that converges vaguely to the Wolpert and Ickstadt (1998) representation of the beta process.\\

This paper is organized as follows. Section 2 discusses metrizing the vague convergence of random measures of the form (\ref{Pn}). It also develops a criterion for which $\mu_n$ converges vaguely to $\mu$.   In Section 3, a finite sum approximation of the beta process is derived.  In Section 4, an example comparing the  performance of the new approximation to other  existing approximations is presented. Section 5 ends with a brief summary of the  results.

\section {Metrizing Vague Convergence of Random Measures}
\label{main}
The material developed in this section  can be seen as a convenient adaptation of the work  of Grandell (1977).  Let $\mu$ be a measure on $\mathbb{R}$ such that its distribution function $\mu(t)=\mu\left((-\infty,t]\right),$ $t \in \mathbb{R},$ is finite for finite $t.$ The same notation will be used for the measure and its distribution function. We assume that $\mu(-\infty)=0.$ The $\mu-$measure of the interval $(a,b]$ is denoted by $\mu\left((a,b]\right)$ or in terms of the distribution function by $\mu(b)-\mu(a).$ The set of measures $\mu$ on $\mathbb{R}$ such that $\mu(-\infty)=0$  and $\mu(t)<\infty$ for all $t \in \mathbb{R}$ is denoted by $\mathcal{M}.$ We shall now define vague convergence on $\mathcal{M}.$

\begin{definition} Let $\mu, \mu_1,\mu_2,\ldots \in \mathcal{M},$ be given. We say that $\mu_n$ converges vaguely to $\mu$ and write that $\mu_n \overset{v} \to\mu$ if  $\mu_n(t) \to\mu(t)$ for all $t \in \mathbb{R}$ such that $\mu$ is continuous at $t.$
\end{definition}
Let $C_K^{+}(\mathbb{R})$  be the set of all  nonnegative continuous real valued functions with compact support defined on $\mathbb{R}.$ In $C_K^{+}(\mathbb{R})$ all functions are bounded and for each $f \in C_K^{+}(\mathbb{R})$ there exists to a number $x_f$ such that $f(x)=0$ for all $x>x_f.$ For the proof of the next theorem consult, for example,  Grandell (1977).

\begin{theorem} \label{thm2} Let  $\mu, \mu_1,\mu_2,\ldots \in \mathcal{M}$ be given. Then, as $n \to \infty$, $\mu_n\overset{v}\to\mu$ if and only if
$$\int_{-\infty}^{\infty} f(x)\mu_n(dx) \to\int_{-\infty}^{\infty} f(x)\mu(dx)$$
for all  $f \in C_K^{+}(\mathbb{R}).$
\end{theorem}

To define a metric, corresponding to vague convergence, in the set $\mathcal{M},$ we first need to define a metric in the following set: $$\mathcal{M}^{(k)}=\left\{\mu \in \mathcal{M}: \mu(x)=\mu(\theta_{(k)}) \text { for } x\ge \theta_{(k)})\right\},$$
where $\theta_{(k)}=\max_{1\le i \le k} \theta_i$. Here, $\theta_i \overset {\text{i.i.d}}\sim \Pi$, where $\Pi$ is  a fixed  continuous probability measure on  $\mathbb{R}.$ In other words, the set $\mathcal{M}^{(k)}$ is the set of those measures in  $\mathcal{M}$ with total mass on the interval $(-\infty,\theta_{(k)}].$ It is worth mentioning that the set $\mathcal{M}^{(k)}$ considered in this section is  different from  that defined in Grandell (1977). Grandell (1977) used  the following set: $$\mathcal{M}^{(k)}_{\text{Grandell}}=\left\{\mu \in \mathcal{M}: \mu(x)=\mu(k) \text { for } x\ge{k}\right\}.$$
Clearly, the random measure $\mu_k$ defined in (\ref{Pn}) puts all its mass in $(-\infty, \theta(k)]$. Thus,  $\mu_k$ belongs to $\mathcal{M}^{(k)}$ but not to $\mathcal{M}^{(k)}_{\text{Grandell}}$.

For $\mu_1$ and $\mu_2 \in \mathcal{M}^{(k)},$ the L\'evy metric  is defined by
 \begin{eqnarray*}
 d_L(\mu_1,\mu_2)=\inf\{h\ge 0: \mu_1(x-h)-h \le \mu_2(x) \le \mu_1(x+h)+h, \forall x \in \mathbb{R} \}.
  \end{eqnarray*}
The next Lemma deals with some properties of $d_L.$ The proof is very similar to the proof of Lemma 1 of Grandell (1977). Thus, the proof is omitted.

\begin{lemma}$d_L$ is a metric in $\mathcal{M}^{(k)},$ i.e. for all $\mu_1,\mu_2,\mu_3 \in \mathcal{M}^{(k)}$ we have
\begin{enumerate}
   \item $d_L(\mu_1,\mu_2)=0$ if and only if $\mu_1=\mu_2.$
   \item $d_L(\mu_1,\mu_2)=d_L(\mu_2,\mu_1).$
    \item $d_L(\mu_1,\mu_3)\le d_L(\mu_1,\mu_2)+d_L(\mu_2,\mu_3)$
    \end{enumerate}
\end{lemma}

The proof of the next lemma follows by imitating the proof of Lemma 2 of Grandell (1977) with $n$ and $k$ are replaced by $\theta_{(k)}$ and $n$, respectively.

\begin{lemma} \label{levy} $d_L$ metrizes vague convergence in $\mathcal{M}^{(k)},$ i.e. $\mu_n \overset{v}\to\mu$ (as $n\to\infty$) if and only if $d_L(\mu_n,\mu) \to 0$ (as $n\to\infty$) for $\mu,\mu_1, \mu_2,\ldots  \in\mathcal{M}^{(k)}.$
\end{lemma}

In Lemma \ref{levy}, we have shown that the L\'evy metric metrizes vague convergence in $\mathcal{M}^{(k)}.$ We will use this to prove a similar result in $\mathcal{M}$. As in Grandell (1977), the idea is to associate to each $\mu \in \mathcal{M}$ a vector $\left(\mu^{(1)},\mu^{(2)},\ldots\right)$, where $\mu^{(k)} \in \mathcal{M}^{(k)}$ and where componentwise convergence is equivalent to convergence. To do this, choose $f_1,f_2,\cdots \in C_K^{+}(\mathbb{R})$ such that
\begin{equation*}
  f_k(t)= \left\{
  \begin{array}{cc}
    1   & t<\theta_{(k-1)} \\
    \theta_{(k)}-t  & \theta_{(k-1)} \le t <\theta_{(k)} \\
    0    & t \ge \theta_{(k)}
  \end{array}
  \right.
\end{equation*}
Define  $\mu^{(k)}(t)=\int_{-\infty}^t f_k(x)\mu(dx).$ Clearly, the mapping $\mu \curvearrowright \left(\mu_{1}^{(k)},\mu_{2}^{(k)},\cdots\right)$ is one to one. Define, for $\mu_1$ and $\mu_2 \in \mathcal{M},$
$$d\left(\mu_1,\mu_2\right)=\sum_{k=1}^\infty \frac {d_L\left(\mu_1^{(k)},\mu_2^{(k)}\right)} {2^k\left(1+d_L\left(\mu_1^{(k)},\mu_2^{(k)}\right)\right)}.$$

\begin{theorem} $d$ is a metric on $\mathcal{M},$ i.e.
\begin{enumerate}
\item $d(\mu_1,\mu_1)=0$ if and only if $\mu=\mu_2$
\item $d(\mu_1,\mu_2)=d(\mu_2,\mu_1).$
\item $d(\mu_1,\mu_3)\le d(\mu_1,\mu_2)+d(\mu_2,\mu_3)$ for all $\mu_1,\mu_2,\mu_3 \in \mathcal{M}.$
\end{enumerate}
\end{theorem}
\proof (i) $d(\mu_1,\mu_1)=0$ if and only if $d_L(\mu_1^{(k)},\mu_2^{(k)})$ if and only if  $\mu_1^{(k)}=\mu_2^{(k)}$ if and only if $\mu_1=\mu_2.$\\

(ii) $d(\mu_1,\mu_2)=d(\mu_2,\mu_1)$ follows directly since $d_L(\mu_1^{(k)},\mu_2^{(k)})=d_L(\mu_2^{(k)},\mu_1^{(k)}).$

(iii) Since   $d_L(\mu_1^{(k)},\mu_3^{(k)})\le d_L(\mu_1^{(k)},\mu_2^{(k)})+d_L(\mu_2^{(k)},\mu_3^{(k)})$ for all $\mu_1^{(k)},\mu_2^{(k)},\mu_3^{(k)} \in \mathcal{M}^{(k)},$ we have
 \begin{eqnarray*}
 \frac {d_L\left(\mu_1^{(k)},\mu_3^{(k)}\right)} {1+d_L\left(\mu_1^{(k)},\mu_3^{(k)}\right)} &=&1-\frac {1} {1+d_L\left(\mu_1^{(k)},\mu_3^{(k)}\right)}\\
&\le& 1-\frac {1} {1+d_L\left(\mu_1^{(k)},\mu_2^{(k)}\right)+d_L\left(\mu_2^{(k)},\mu_3^{(k)}\right)}\\
 &=& \frac {d_L\left(\mu_1^{(k)},\mu_2^{(k)}\right)+d_L\left(\mu_2^{(k)},\mu_3^{(k)}\right)} {1+d_L\left(\mu_1^{(k)},\mu_2^{(k)}\right)+d_L\left(\mu_2^{(k)},\mu_3^{(k)}\right)}\\
 &=&\frac {d_L\left(\mu_1^{(k)},\mu_2^{(k)}\right)} {1+d_L\left(\mu_1^{(k)},\mu_2^{(k)}\right)+d_L\left(\mu_2^{(k)},\mu_3^{(k)}\right)}\\&&+\frac {d_L\left(\mu_2^{(k)},\mu_3^{(k)}\right)} {1+d_L\left(\mu_1^{(k)},\mu_2^{(k)}\right)+d_L\left(\mu_2^{(k)},\mu_3^{(k)}\right)}\\
 &\le& \frac {d_L\left(\mu_1^{(k)},\mu_2^{(k)}\right)} {1+d_L\left(\mu_1^{(k)},\mu_2^{(k)}\right)}+\frac {d_L\left(\mu_2^{(k)},\mu_3{(k)}\right)} {1+d_L\left(\mu_2^{(k)},\mu_3^{(k)}\right)},
\end{eqnarray*}
Thus,
 \begin{eqnarray*}
 d\left(\mu_1,\mu_3\right)&=&\sum_{k=1}^\infty \frac {d_L\left(\mu_1^{(k)},\mu_3^{(k)}\right)} {2^k\left(1+d_L\left(\mu_1^{(k)},\mu_3^{(k)}\right)\right)}\\
&\le& \sum_{k=1}^\infty \frac {d_L\left(\mu_1^{(k)},\mu_2^{(k)}\right)} {2^k\left(1+d_L\left(\mu_1^{(k)},\mu_2^{(k)}\right)\right)}+\sum_{k=1}^\infty \frac {d_L\left(\mu_2^{(k)},\mu_3^{(k)}\right)} {2^k\left(1+d_L\left(\mu_2^{(k)},\mu_3^{(k)}\right)\right)} \\
 &=&d\left(\mu_1,\mu_2\right)+d\left(\mu_2,\mu_3\right).
\end{eqnarray*}
 This completes the proof of the lemma.
 \endproof

\begin{theorem}
$d$ metrizes vague convergence in $\mathcal{M}.$
\end{theorem}

\proof Let $\mu,\mu_1,\mu_2, \ldots \in \mathcal{M}$  be given. It follows from the definition of $d$ that $d(\mu_n,\mu) \to0$ if and only if $d_L(\mu_n^{(k)},\mu^{(k)})\to0$ (as $n\to\infty$) for all $k$.  By Lemma \ref{levy}, this holds if and only if $\mu_n^{(k)} \overset{v} \to\mu^{(k)}$ for all $k$. Thus, it is enough to prove that $\mu_n \overset{v} \to\mu$ if and only if $\mu_n^{(k)} \overset{v} \to\mu^{(k)}$ for all $k.$

By Theorem \ref{thm2},  $\mu_n \overset{v} \to\mu$ if and only if $\int f(x)\mu_n(dx)\to\int f(x)\mu(dx)$ for all $f \in C_K^{+}(\mathbb{R}).$ Since $ff_k \in C_K^{+}(\mathbb{R})$ for all $k$ and all $f\in C_K^{+}(\mathbb{R}),$ it follows that $\mu_n \overset{v}\to\mu$ implies that $\mu_n^{(k)} \overset{v}\to\mu^{(k)}$ for all $k.$

Conversely, if $\mu_n^{(k)} \overset{v}\to\mu^{(k)}$ for all $k$ we can for each $f \in C_K^{+}(\mathbb{R})$ choose $k$ so that $ff_k=f.$ Thus, we have
\begin{eqnarray}
\nonumber \int f(x)\mu_n(dx)&=& \int f(x)f_k(x)\mu_n(dx)= \int f(x)\mu_n^{(k)}(dx)\\
\nonumber  &\to& \int f(x)\mu^{(k)}(dx)=\int f(x)\mu(dx).T
\end{eqnarray}
Thus, $\mu_n \overset{v}\to\mu.$
\endproof

The proof of Theorem 3 reveals the following interesting result.

\begin{corollary}\label{corollary1} Let $\mu,\mu_1,\mu_2, \ldots \in \mathcal{M}$ and $\mu^{(k)},\mu_1^{(k)},\mu_2^{(k)}, \ldots \in \mathcal{M}^{(k)}$. Then $\mu_n \overset{v}\to\mu$ (as $n  \to\infty$) if and only if $\mu_{n}^{(k)} \overset{v} \to\mu^{(k)}$ (as $n  \to\infty$) for all $k$ fixed.
\end{corollary}

\section{Applications on Bayesian Nonparametric Priors}
\label{applications}
An interesting application of  Corollary \ref{corollary1}   comprises deriving a finite sum representation that converge vaguely  to the Wolpert and Ickstadt (1998) representation of the beta process  (Hjort, 1990).  Finite sum approximations for the Dirichlet process, beta-Stacy process, normalized inverse-Gaussian process and two-parameter Poisson-Dirichlet process were derived, respectively,  in Zarepour and Al-Labadi (2012) and Al-Labadi  and Zarepour (2013a,b; 2014a,b).

There are two common techniques to writing a series representation for any L\'evy process having no Gaussian component.
The first comes from Ferguson and  Klass (1972). The second technique is from Wolpert and Ickstadt (1998).  A brief discussion  of the two  methods is described in Appendix A of this paper. An interesting comparison between the two representations from the computational point of view was addressed in  Al-Labadi  and Zarepour  (2013a). Here, it is pointed out that the representation of Wolpert and Ickstadt is more appropriate for dealing with  nonhomogeneous processes (i.e., the L\'evy measure in (\ref{B5}) depends on $t$). Conversely, for homogeneous processes (i.e., the L\'evy measure is independent of $t$), the two approaches are equivalent.

Let $A_0$ be a continuous cumulative hazard function and $c(t)_{t\ge 0}$
be a piecewise continuous, nonnegative function. Following Hjort (1990), the beta process  $A$, written $A \sim BP \left(c(\cdot),A_0(\cdot)\right)$, is the completely random measure with L\'evy measure
\begin{flalign}
L_t(ds)& =\left[\int_{0}^{t}c(z)s^{-1}(1-s)^{c(z)-1}dA_{0}(z)\right]ds, \ \ \text{for } t \ge 0,\ \ 0<s<1. \label{B5}&
\end{flalign}

By (\ref{ferg}) and (\ref{WI}),  since no closed form for the inverse of the L\'evy measure (\ref{B5}) exists, the simulation of the  beta process based on series representations is very complex and may be difficult to apply in practice for many users. The next theorem outlines a remedy to this problem. Note that, when $c(t)=c$ for all $t$ (i.e. the homogenous case), Al-Labadi and Zarepour (2015) derived a finite sum approximation and showed that it converges almost surely to the representation of Ferguson and  Klass (1972) of the beta process. More details about interesting properties of the beta process when $c(t)=c$ for all $t$ are discussed Al-Labadi and Abdelrazeq (2016).

\begin{theorem} \label{Beta}
Let $(\theta_i)_{i \geq 1}$ be i.i.d. random variables with common distribution $\Pi$ and $\Gamma_i=E_1+\cdots+E_i,$ where $\left(E_i\right)_{i\ge1}$ are i.i.d. with exponential distribution with mean 1, independent of $(\theta_i)_{i\ge1}.$ Let $A \sim {BP}\left(c(\cdot),A_0(\cdot)\right)$ on $[0,t_0]$,  where $t_0>0$ is fixed. We assume that $A_0$ is continuous with $A_0(t_0)< \infty$. Let $\Pi(dz)=\eta(dz)/A_0(t_0),$ where $\eta([0,t])=A_0(t)$.
\begin{eqnarray}
\label{Levy1}  L_{n,\theta}(x)&=&\frac{\Gamma\left(c(\theta)\right)}{\Gamma({c(\theta)}/{n})\Gamma\left(c(\theta)-{c(\theta)}/{n}\right)}\int_{x}^1
s^{{c(\theta)}/{n}-1}\left(1-s\right)^{c(\theta)\left(1-{1}/{n}\right)-1}ds. \ \ \ \  \  \
\end{eqnarray}
and
\begin{eqnarray*}
  M_{z}(x)=A_0(t_0)c(z)\int_{x}^1 s^{-1}(1-s)^{c(z)-1}ds.
\end{eqnarray*}

Then, as $n\to\infty$,
\begin{equation}
A_{n}(t)=\sum_{i=1}^{n} {{L^{-1}_{n,\theta_i}\left(\frac{\Gamma_i}{ A_0(t_0)n}\right)}\delta_{\theta_i}} \overset{v}\rightarrow
A(t)=\sum_{i=1}^{\infty} {M^{-1}_{\theta_i}\left(\Gamma_i\right)\delta_{\theta_i}}. \label{beta}
\end{equation}
\end{theorem}

\proof First we show that, for any  $x \in (0,1)$,
\begin{equation}
\label{Levy3} nA_0(t_0)L_{n,\theta}(x)\to M_{\theta}(x).
\end{equation}
Note that, for any $x > 0$,
\begin{equation*}
\Gamma(x)=\frac {\Gamma(x+1)}{x}.
\end{equation*}
With $x=c(\theta)/n$ we obtain
\begin{equation*}
\frac {n}{\Gamma(c(\theta)/n)}=\frac {c(\theta)}{\Gamma(c(\theta)/n+1)}.
\end{equation*}
Since
$\Gamma(x)$ is a continuous function,
$$\frac{n}{\Gamma(c(\theta)/n)} \times \frac{\Gamma\left(c(\theta)\right)}{\Gamma\left(c(\theta)-{c(\theta)}/{n}\right)}\to c(\theta).$$
Clearly, the integrand in the right hand side of (\ref{Levy1}) converges to $s^{-1}(1-s)^{c(\theta)-1}.$ To apply the dominated convergence theorem, we need to show that  this integrand  is dominated   by an integrable function. Since $x<s<1$, we have $s^{-1}<x^{-1}$ and $s^{{c(\theta)}/{n}}<1.$ This implies that $s^{ {c(\theta)}/{n}-1}<x^{-1}.$ Therefore, the integrand is bounded above by the integrable function $x^{-1}\left(1-s\right)^{c(\theta)
\left(1-{1}/{n}\right)-1}$. Thus, by the dominated convergence theorem, we get  (\ref{Levy3}). Since that the left hand side of (\ref{Levy3}) is a sequence of a continuous monotone functions converging to a monotone function for every $x>0$. This is equivalent to the convergence of their inverse function to the inverse function of the right hand side (de Haan \& Ferreira, 2006, page 5). Thus, as $n\rightarrow \infty$,
$$L^{-1}_{n,\theta_i}\left(\frac{\Gamma_i}{n A_0(t_0)}\right)  \overset{v}\rightarrow M^{-1}_{\theta_i}(\Gamma_i),$$
To complete the proof of the theorem, we apply  Corollary \ref{corollary1} with
\begin{equation}
\mu_{n}^{(k)}=A_{n}^{(k)}=\sum_{i=1}^{k} L^{-1}_{n,\theta_i}\left(\frac{\Gamma_i}{A_0(t_0)n}\right)\delta_{\theta_{(i)}} \nonumber
\end{equation}
and
\begin{equation}
\mu^{(k)}=A^{(k)}=\sum_{i=1}^{k} {M^{-1}_{\theta_i}}\left(\Gamma_i\right)\delta_{\theta_{(i)}}. \nonumber
\end{equation}
Clearly, both $A_{n}^{(k)}$  and $A^{(k)}$ belong to $\mathcal{M}^{(k)}$. Since, for all $k$ fixed,
\begin{equation}
A_{n}^{(k)} \overset {v}\to A^{(k)}, \nonumber
\end{equation}
as $n  \to\infty$, we get (\ref{beta}).  This completes the proof of the theorem.
\endproof

Note that, $L^{-1}_{n,\theta_i}\left(p\right)$ is the $1-p$-th quantile of ${beta}\left(c(\theta_i)/n,c(\theta_i)(1-1/n)\right)$ distribution. This provides the following algorithm.
\begin{enumerate}
\item Fix a relatively large positive integer $n$.
\item For $i=1,\ldots,n$, generate  $\theta_i\overset{\text{i.i.d.}}\sim \Pi$, where $\Pi(dz)=\eta(dz)/A_0(t_0)$ and $\eta([0,t])=A_0(t)$.
\item For  $i=1,\ldots,n+1$, generate $E_i\overset{\text{i.i.d.}}\sim exponential(1)$   such that $\left(E_i\right)_{1\le i \le n+1}$ and $\left(\theta_i\right)_{1\le i \le n}$ are independent.
\item For $i=1, \ldots,n+1,$ compute $\Gamma_i=E_1+\cdots+E_i.$
\item For $i=1, \ldots,n,$ compute $L^{-1}_{n,\theta_i}\left(\Gamma_i/\left(A_0(t_0)n\right)\right)$.
\item Use (\ref{beta}) to obtain an approximate value of $A \sim BP(c(\cdot), A_0(\cdot))$.
 \end{enumerate}

Note that, it is possible to extend Theorem \ref{Beta} to derive an approximation of the  beta-Dirichlet process   (Kim,  James and Weibbach, 2012), a nonparametric prior for the cumulative intensity functions of a Markov process.  Specifically, as $n\to \infty$,
\begin{equation}
\nonumber B_{n}(t)=\sum_{i=1}^{n} V(\theta_i){{L^{-1}_{n,\theta_i}\left(\frac{\Gamma_i}{ A_0(t_0)n}\right)}\delta_{\theta_i}} \overset{v}\rightarrow
B(t)=\sum_{i=1}^{\infty} V(\theta_i){M^{-1}_{\theta_i}\left(\Gamma_i\right)\delta_{\theta_i}},
\end{equation}
where $L^{-1}_{n,\theta_i}$, $M^{-1}_{\theta_i}\left(\Gamma_i\right)$ are  defined as in Theorem  \ref{Beta} and $V(\theta_i)$ are independent Dirichlet of random vectors with parameters $\gamma_1(\theta_i), \ldots, \gamma_n(\theta_i)$. Here $B$ is the beta-Dirichlet process with parameters $(A_0,c, \gamma_1,\gamma_1, \ldots, \gamma_n)$. We refer the reader to the paper of Kim,  James and Weibbach (2012) for the details.

\section{Empirical Results: Comparison to Other Methods}

Sampling from the beta process plays a central role in many applications. We refer the reader to the work of  Paisley and Carin (2009) and Broderick, Jordan, and Pitman (2012). It is also  required to simulate the beta-Dirichlet process (Kim,  James and Weibbach, 2012). Several algorithms  to sample from the beta process   exist in the literature. In this section, we compare the new approximation of the beta process with  the algorithm of Ferguson and  Klass  (1972),  the algorithm of Damien, Laud, and Smith (1995), the algorithm of Wolpert and Ickstadt (1998), the algorithm of Lee and Kim  (2004) and  the algorithm of Lee (2007). A summary of these algorithms is given in Appendix A.

In order to make comparisons between the  preceding algorithms, we  use equivalent settings for the parameters characterizing these algorithms (see Table 1).   We consider the beta process with $c(t)=2e^{-t}$ and   $A_0(t) = t$, where $t \in [0,1]$.  We  compute  the absolute maximum difference  between an approximate sample mean and the exact mean. See also Lee and Kim (2004) and Lee (2007) for similar comparisons.  The exact mean of $A(t)$ in this example is $A_0(t)=t$; see Hjort (1990). We refer to this statistic by the maximum mean error. Specifically,
\begin{eqnarray*}
\text{maximum mean error}=\max_{t}\left|E\left[A_n(t)\right]-E\left[A(t)\right]\right|&=& \max_{t}\left|E\left[A_n(t)\right]-t\right|,
\end{eqnarray*}
where $t=0.1,0.2,\ldots,0.9,1.0$ and $A_n$ is an approximation of $A\sim {BP}(c(t)=2e^{-t},A_0(t)=t).$
Note that $E\left[A_n(x)\right]$ is approximated  by obtaining  the mean  at $t=0.1,0.2,\ldots,0.9,1.0$ of $3000$ i.i.d. sample paths from the approximated process $A_n$. Similarly, we compute the maximum standard deviation  error between an approximate sample standard deviation (s.d.) and the exact standard deviation. The exact standard deviation of $A(t)$  is $\sqrt{t/3}$; see Hjort (1990). Thus,
\begin{eqnarray*}
\text{maximum s.d. error}&=&\max_{t}\left|s.d\left[A_n(t)\right]-s.d.\left[A(t)\right]\right|\\
&=& \max_{t}\left|s.d.\left[A_n(t)\right]-\sqrt{t/3}\right|.
\end{eqnarray*}
Table 1 depicts values of the maximum mean  error, the maximum standard deviation  error, and the corresponding computational time. Simulating the algorithm of Ferguson and  Klass  (1972) and the algorithm of Wolpert and Ickstadt (1998) is performed through relatively complex  numerical methods, which are not appropriate for many users and  time consuming (See Table 1). The R function  ``uniroot" is used to implement these two algorithms. The computational time is computed by applying the R function ``System.Time". As seen in Table 1, the new algorithm has the smallest mean and standard deviation errors. Furthermore, it has a very reasonable computation time.

\begin{table}[htbp]
\caption{This table reports the maximum mean error, the maximum standard deviation error, and the corresponding computation time. Here, FK, DSL, WI and LK stand for  the algorithm of Ferguson and Klass (1972), the algorithm of Damien, Laud, and Smith (1995), the algorithm of Wolpert and Ickstadt (1998), and the algorithm of Lee and Kim (2004), respectively.}
\begin{center}
\begin{tabular}{lllll}
\hline
\hline
\multicolumn{0}{c} {Algorithm} &\multicolumn{0}{c} {Parameters} & \multicolumn{0}{c} {max. mean error}  &
\multicolumn{0}{c} {max. s.d. error}& \multicolumn{0}{c} {Time}\\
\hline
KL &$n=200$ & 0.0192& 0.1047&657.53\\
\\
DSL &$m=n=200$ & 0.0167& 0.0145&90.53\\
\\
WI &$M=200$   & 0.0167 & 0.0884&406.10 \\
\\
LK  &$\epsilon=0.01$   &   0.0217&0.0239& 0.29\\
\\
Lee  &$n=200, \epsilon=0.05$ & 0.0125&0.0522&1.25\\
\\
New &$n=200$   & 0.0069&0.0089&6.28\\
\\
\hline
\end{tabular}
\end{center}
\label{table4.3}
\end{table}

\section {Conclusions}

The vague convergence of random measures of the form  (\ref{Pn}) has been studied in this paper. An interesting application of the derived results includes deriving a finite sum representation that converges vaguely  to the representation of Wolpert and Ickstadt (1998)  of the beta process. This representation gives a simple yet efficient approach to approximate the beta process.   We believe that the comprehensive study of metrizing random measures as in (\ref{Pn}) and its strong association to various  Bayesian nonparametric priors will add further useful tools to the Bayesian nonparametric toolbox.

\section{Acknowledgments}
 Research of the author is supported by the \textbf{Natural Sciences and Engineering Research Council of Canada (NSERC)}.


\appendix

\section {Other Sampling Algorithms} Below is a brief discussion  of the algorithms considered in Section 4 of the present paper. We refer the reader to the original papers for more details. Let $A \sim {BP}\left(c(\cdot),A_0(\cdot)\right)$ on $[0,t_0]$,  where $t_0>0$ is fixed.  We assume that $A_{0}(t_0)<\infty$.

\smallskip

\noindent \textbf{ $\bullet$ Ferguson-Klass Algorithm:} The  steps of the algorithm of Ferguson and Klass (1972) are:

\begin{enumerate}
\item Let $\Gamma_i=E_1+\cdots+E_i,$ where $(E_i)_{i\ge 1}$ are i.i.d. random variables with exponential distribution of mean 1.

\item For each $i \ge 1,$  let $J_i$ be the solution of  $\Gamma_i=L_{t_0}\left(J_i\right),$ where $L_{t_0}\left(x\right)=L_{t_0}\left([x,1)\right),$ $x>0$
and the measure $L_t$ is given by (\ref{B5}).
\item Generate i.i.d. random variables $(U_i)_{i\ge 1}$  from the uniform distribution on $[0,1],$ independent of $(E_i)_{i\ge 1}.$

\item For $i \ge 1,$  let $\theta_i$ be the solution of $U_i=n_{\theta_i}(J_i)$ in $[0,t_0]$ where
\begin{equation}
\nonumber n_t(s)=\frac{\int_{0}^t c(z)s^{-1}(1-s)^{c(z)}dA_{0}(z)}{\int_{0}^{t_0} c(z)s^{-1}(1-s)^{c(z)}dA_{0}(z)}=\frac{\int_{0}^t c(z)(1-s)^{c(z)}dA_{0}(z)}
{\int_{0}^{t_0} c(z)(1-s)^{c(z)}dA_{0}(z)}.
\end{equation}
\end{enumerate}
\smallskip
The process
\begin{equation}
A(t)=\sum_{i=1}^\infty J_i  I(\theta_i \le t)=\sum_{i=1}^\infty L_{t_0}^{-1}(\Gamma_i)  I(\theta_i \le t) \label{ferg}
\end{equation}
is a beta process with parameters $c(\cdot)$ and $A_0(\cdot).$ This series is an infinite series.
In practice, we truncate this series and use the approximation
\begin{equation}
A_n(t)=\sum_{i=1}^n J_i  I(\theta_i \le t)=\sum_{i=1}^n L_{t_0}^{-1}(\Gamma_i)  I(\theta_i \le t). \nonumber
\end{equation}

\noindent \textbf{ $\bullet$ Damien-Laud-Smith Algorithm:} Using  the fact that the distributions of the increments of a nondecreasing L\'evy process
are infinitely divisible, Damien, Laud, and Smith (1995) derived an algorithm to generate approximations for infinitely divisible random
variables and used it to generate  the beta process. First, the time interval $[0,t_0]$ is partitioned into small subintervals with endpoints
$0=\theta_0<\theta_1<\ldots<\theta_m=t_0.$ Let $p_i$ denotes the increment of the process $A$ in  the interval $\Delta_i=(\theta_{i-1},\theta_i],$ i.e. $p_i=A(\theta_i)-A(\theta_{i-1}).$ The steps of the Damien-Laud-Smith algorithm for simulating  the beta process are:
\begin{enumerate}
\item [(1)] Fix a  relatively large positive integer $n$.
\item [(2)] Generate  independent values $z_{ij}$ from $\Pi$, where $\Pi(dz)=\eta(dz)/A_0(t_0)$ and $\eta([0,t])=A_0(t)$, for $j=1,\ldots,n.$
\item [(3)] Generate $x_{ij} \sim {beta} (1,c(z_{ij})),$ for $j=1,\dots,n.$
\item [(4)] Generate $y_{ij}$: $y_{ij}|x_{ij} \sim {Poisson}(\lambda_in^{-1}x_{ij}^{-1})$, for $j=1,\dots,n,$ where $\lambda_i=A_0(\theta_i)-A_0(\theta_{i-1})$.
\item [(5)] Set $p_{i,n}=\sum _{j=1}^n x_{ij}y_{ij}.$ For large $n$, $p_{i,n}$ is an approximation of $p_i$.
\item [(5)] Set $A_{n}=\sum_{i=1}^n p_{i,n} \delta_{\theta_{(i)}}.$
\end{enumerate}
\smallskip
Damien, Laud, and Smith (1995) showed that $A_{n} \overset{d}\to A,$ as $n \to \infty$.

\noindent \textbf{ $\bullet$  Wolpert-Ickstadt Algorithm:} The  steps of the algorithm of  Wolpert and Ickstadt (1998) are:
\begin{enumerate}
\item For $i=1,2,\ldots$, generate  $\theta_i  \overset {i.i.d.} \sim \Pi$, where $\Pi(dz)=\eta(dz)/A_0(t_0)$ and $\eta([0,t])=A_0(t)$.
\item Let $\Gamma_i=E_1+\cdots+E_i,$ where $(E_i)_{i\ge 1}$ are i.i.d. random variables with exponential distribution of
 mean 1, independent of $\left(\theta_i\right)_{i \ge 1}$.

\item Define
 \begin{equation*}
  M_{z}\left(x\right)=\int_{x}^\infty A_0({t_0})c(z)s^{-1}(1-s)^{c(z)-1}ds.
 \end{equation*}

\item For each $i \ge 1,$ solve the equation
\begin{eqnarray*}
  M_{\theta_i}\left(J_i\right)&=&\Gamma_i
\end{eqnarray*}
for $J_i,$ where $\Gamma_i=E_1+\cdots+E_i,$  $(E_i)_{i\ge 1}$ are i.i.d. random variables with exponential distribution of  mean 1 and independent of $\left(\theta_i\right)_{i\ge 1}.$
\item Set
\begin{equation}
 \begin{split}
A(t)&=\sum_{i=1}^\infty M^{-1}_{\theta_i}(\Gamma_i) I\left\{\theta_i\le t\right\}. \label{WI}
 \end{split}
\end{equation}
\end{enumerate}
\smallskip

The process $A$  in (\ref{WI}) is a beta process with parameters $c(\cdot)$ and $A_0(\cdot).$ This series is an infinite series.
In practice, we truncate this series and use the approximation
\begin{equation}
A_n(t)=\sum_{i=1}^n M^{-1}_{\theta_i}(\Gamma_i) I\left\{\theta_i\le t\right\}. \nonumber
\end{equation}

\noindent \textbf{$\bullet$ Lee-Kim Algorithm:}  First the L\'evy measure (\ref{B5}) of the beta process  is approximated by
\begin{equation}
L_{t,\epsilon}(ds)=\left[\int_{0}^t\frac{c(s)}{\epsilon} b(s:\epsilon,c(z))dA_{0}(z)\right]ds, \nonumber
\end{equation}
where
\begin{equation}
b(x:a,b)=\frac{\Gamma(a+b)}{\Gamma(a)\Gamma(b)}x^{a-1}(1-x)^{b-1}, \quad \text{for } 0<x<1, a>0, b>0.\label{B7}
\end{equation}

The steps of the  algorithm of Lee and Kim (2004) for the beta process  are:
\begin{enumerate}
\item [(1)] Fix a  relatively small positive number $\epsilon$.
\item [(2)] Generate the total number of jumps  $n\sim {Poisson} \left(\lambda_\epsilon\right)$, where
$\lambda_\epsilon=L_{t_0,\epsilon}\left((0,1)\right)=\epsilon^{-1}\int_{0}^{t_0}\int_{0}^1 c(z) b(s:\epsilon,c(z))dsdA_{0}(z)=\epsilon^{-1} \int_{0}^{t_0} c(z)dA_{0}(z)$
\item [(3)] Generate the jump times $\theta_1,\ldots,\theta_n$ form the probability density function $dG_\epsilon/\lambda_\epsilon,$ where $dG_\epsilon(z)=  \epsilon^{-1} c(z)  dA_{0}(z)I(0\le z \le {t_0}).$

\item [(4)] Let $\theta_{(1)} \le \ldots \le \theta_{(n)}$ be the corresponding order statistics of $\theta_1, \ldots, \theta_n.$
\item [(5)] Generate the jump sizes $p_1,\ldots,p_n:$ $p_i|\theta_{(i)} \sim {Beta}(\epsilon,c(\theta_{(i)})).$
\item [(6)] Set $A_{\epsilon}=\sum_{i=1}^n p_i \delta_{\theta_{(i)}}.$
\end{enumerate}
\smallskip

Lee and Kim (2004) showed that $A_{\epsilon} \overset{d}\to A,$ as $\epsilon \to 0$.

\smallskip

\noindent \textbf{$\bullet$ Lee  Algorithm:} The steps of the  algorithm of Lee (2007) are:

\begin{enumerate}
\item [(1)] Fix a  relatively large positive integer $n$.
\item [(2)] For $i=1,2,\ldots,n$, generate  $\theta_i  \overset {i.i.d.} \sim \Pi$, where $\Pi(dz)=\eta(dz)/A_0(t_0)$ and $\eta([0,t])=A_0(t)$.
\item [(3)] For $i=1,\dots,n,$ generate $x_i \sim b(s:\epsilon,c(\theta_i)),$ where $b(s:\epsilon,c(\theta_i))$ is defined in (\ref{B7}).
\item [(4)] For $i=1,\dots,n,$ generate $y_i\sim Poisson\left(A_{0}(t_0)b(x_i:1,c(\theta_i))/(n x_i b(x_i:\epsilon,c(\theta_i)))\right).$
\item [(5)] Set $A_{n}=\sum_{i=1}^n x_i y_i \delta_{\theta_i}.$
\end{enumerate}
\smallskip

Lee (2007) proved that, as $n\to \infty,$  $A_{n} \overset{d} \to A$.

\end{document}